\documentclass{amsart}

\usepackage{amssymb}
\usepackage{times}
\usepackage{mathrsfs}

\usepackage[colorlinks,linkcolor={blue},
citecolor={blue},urlcolor={red},]{hyperref}

\theoremstyle{plain}
\newtheorem{theorem}{Theorem}[section]
\theoremstyle{remark}

\theoremstyle{plain}

\newtheorem{lemma}[theorem]{Lemma}

\numberwithin{equation}{section}

\font\Bbbten=msbm10
\font\Bbbseven=msbm7
\font\Bbbfive=msbm5
\newfam\Bbbfam

  \textfont\Bbbfam=\Bbbten
  \scriptfont\Bbbfam=\Bbbseven
  \scriptscriptfont\Bbbfam=\Bbbfive

\def\N{{\mathbb N}}

\newcommand{\e}{\varepsilon}
\newcommand{\D}{\delta}
\newcommand{\s}{^{\ast}}

\newcommand{\n}{\Vert}

\newcommand{\lb}{\langle}
\newcommand{\rb}{\rangle}

\renewcommand{\d}{\delta}

\newcommand{\limn}{\lim_{n\to\infty}}
\newcommand{\dps}{\displaystyle}

\begin{document}
\title
{Separated sequences in uniformly convex Banach spaces}

\author{J.M.A.M. van Neerven}
\address{Delft Institute of Applied Mathematics\\
Technical University of Delft \\ P.O. Box 5031\\ 2600 GA Delft\\The Netherlands}
\email{J.vanNeerven@math.tudelft.nl}



\begin{abstract}
We give a characterization of uniformly convex Banach spaces in
terms of a uniform version of the Kadec-Klee property. As an
application we prove that if $(x_n)$ is a bounded sequence
in a uniformly convex Banach space $X$ which is $\e$-separated for some $0<\e\le 2$,
then for all norm one vectors $x\in X$ there
exists a subsequence $(x_{n_j})$ of $(x_n)$ such that
$$\dps \inf_{j\not=k}\ \bigl\n x-(x_{n_j} - x_{n_k})
\bigr\n \ge 1+\D_X(\tfrac23\e),$$ where $\D_X$ is
the modulus of convexity of $X$. From this we deduce that
the unit sphere of every infinite-dimensional uniformly convex Banach space 
contains a $(1+\tfrac12\D_X(\tfrac23))$-separated 
sequence.
\end{abstract}

\thanks{This work was supported by
a `VIDI subsidie' in the `Vernieuwingsimpuls' programme
of the Netherlands Organization for Scientific Research (NWO)
and by the Research Training Network HPRN-CT-2002-00281}
\maketitle

\section{Introduction and statement of the results}

In this note we obtain a characterization of uniformly convex
Banach spaces in terms of a uniform version of the Kadec-Klee property.
Recall that a Banach space $X$ with unit sphere $S_X$ is called
{\em uniformly convex} \cite{De} if for all $0<\e\le 2$ we have
$\D_X(\e)>0$, where
$$\e\mapsto \D_X(\e):= \inf\big\{ 1- \tfrac12\n x+y\n: \ x,y\in S_X, \
\n x-y\n\ge\e\big\}
$$
denotes the modulus of convexity of $X$.
For general properties of this function we refer to \cite{De,Di}.
Before stating the main abstract result of this paper we formulate two applications.
The first concerns the differences $x_j-x_k$ of a uniformly
separated bounded sequence $(x_n)$. To motivate the result let us
first consider an arbitrary bounded sequence $(x_n)$ in a Banach
space $X$. It is easy to see that for all $x\in S_X$ there exists
a subsequence $(x_{n_j})$ of $(x_n)$ such that
$$ \liminf_{j,k\to\infty}\ \bigl\n x-(x_{n_{j}} -x_{n_k})
\bigr\n \ge 1.$$
Indeed, let $x\in S_X$ be arbitrary and
choose $x\s\in S_{X\s}$ such that $|\lb x, x\s\rb|=1$. If the 
subsequence $(x_{n_j})$ is chosen in such
a way that the scalar sequence $(\lb x_{n_j}, x\s\rb)$ is convergent, then
$$ \liminf_{j,k\to\infty}\ \bigl\n x-(x_{n_{j}} -x_{n_k})
\bigr\n \ge  \liminf_{j,k\to\infty}\ \bigl|
\lb x-(x_{n_{j}} -x_{n_k}), x\s\rb\bigr| = |\lb x,x\s\rb| =
1.$$ Easy examples show that for $X=C[0,1]$ the value of the constant
$1$ is the best possible, even if the sequence $(x_n)$ is assumed
to be $\e$-separated for some $\e>0$, by which we mean
that $\n x_j - x_k\n\ge \e$ for all $j\not=k$. For uniformly convex spaces $X$
we obtain an improved constant in terms of the modulus of convexity $\D_X$:

\begin{theorem}\label{thm:main1}
Let $X$ be a uniformly convex Banach space and 
let $(x_n)$ be a bounded 
sequence in $X$ which is $\e$-separated for some $0<\e\le 2$. Then for all $x\in S_X$ there exists a 
subsequence $(x_{n_j})$ of $(x_n)$ satisfying
$$\inf_{j\not=k}\ \bigl\n x-(x_{n_j} - x_{n_k})
\bigr\n \ge 1+ \D_X(\tfrac23 \e).$$
\end{theorem}

A celebrated result of Elton and Odell \cite{EO} asserts
that the unit sphere of  every infinite-dimensional Banach space
$X$ contains a $(1+\mu)$-separated sequence for some $\mu>0$.
It was subsequently shown by
Kryczka and Prus \cite{KP} that the unit sphere of every
infinite-dimensional nonreflexive Banach space contains a
$\sqrt[5]{4}$-separated sequence. 
For uniformly convex spaces $X$ we use Theorem \ref{thm:main1} 
to deduce a lower bound for the separation constant 
in terms of the modulus of convexity $\D_X$:

\begin{theorem}\label{thm:main3}
The unit sphere of every infinite-dimensional uniformly convex
Banach space $X$ contains a
$(1+\tfrac12\D_X(\tfrac23))$-separated sequence.
\end{theorem}

Since uniformly convex spaces are reflexive, this result does
not overlap with the result of Kryczka and Prus.
Theorem \ref{thm:main3} 
provides an affirmative answer, for the class of uniformly convex
spaces, to a question of Diestel \cite[page 254]{Di}.

In $X=l^p$, the sequence of unit vectors is $2^{\frac1p}$-separated.
On the other hand it was shown by Clarkson \cite{Cl} and Hanner \cite{Ha} that 
$l^p$ is uniformly convex for $p\in (1,\infty)$ 
with modulus of convexity given by 
$$
\D_{l^p}(\e) = 1 - \big(1-(\tfrac12\e)^p\big)^\frac1p
$$
if $p\in [2,\infty)$, and by the equation 
$$
 \big|1-\D_{l^p}(\e)+\tfrac12\e\big|^p + \big|1-\D_{l^p}(\e)-\tfrac12\e\big|^p =2
$$ 
if $p\in (1,2]$.
Thus for the spaces $l^p$,
Theorem \ref{thm:main3} does not give the best possible separation constant. 
This raises the question whether 
Theorem \ref{thm:main3} can be further improved.

\medskip

Theorems \ref{thm:main1} and \ref{thm:main3} 
are obtained as consequences of the following characterization of
uniformly convex spaces and a quantitative refinement stated in the next
section.

\begin{theorem}\label{thm2}
For a Banach space $X$ the following assertions are equivalent:
\begin{enumerate}
\item $X$ is uniformly convex;

\item For all $\e>0$ there exists a constant
$\d>0$ such that for all $x\in S_X$, all $x'\in X$,
and all linear contractions $T$ from
$X$ into some Banach space $Y$ satisfying
\begin{enumerate}
\item[]
\begin{enumerate}
\item[(i)] $\dps \big|1-\n x'\n \big| < \d$, 
\item[(ii)] $\dps \n T x\n  > 1- \d$ 
\item[(iii)] $\dps \n Tx - Tx'\n < \d,$
\end{enumerate}
\end{enumerate}
we have $ \n x-x'\n < \e$;

\item For all $\e>0$ there exists a constant $\d>0$ such that
for all $x\in S_X$, all $x'\in S_X$, and all
$x\s\in S_{X\s}$ satisfying
\begin{enumerate}
\item[]
\begin{enumerate}
\item[(iv)] $\dps |\lb x, x\s\rb| > 1-\d,$
\item[(v)] $\dps |\lb x-x', x\s \rb| < \d$,
\end{enumerate}
\end{enumerate}
we have $ \n x-x'\n < \e.$
\end{enumerate}
\end{theorem}

Condition  (3), when reformulated in terms of sequences, may be
regarded as a uniform version of
the Kadec-Klee property. Recall that $X$ is said to have the {\em
Kadec-Klee property} if for all $x\in S_X$ and all sequences $(x_n)\subseteq
S_X$ such that $\limn \lb x-x_n, x\s \rb= 0$ for all $x\s\in X\s$
we have $\limn \n x - x_n\n = 0$. Every
uniformly convex space has the Kadec-Klee property; this fact is
also contained as a special case in Theorem \ref{thm2}.

\section{proofs}

\begin{proof}[Proof of Theorem \ref{thm2}.]

(1)$\Rightarrow$(2): \ It is enough to prove this implication for $0<\e\le 2$.
For such an $\e$, we will show that (2) holds for any $\d>0$ satisfying
\begin{equation}\label{delta} \d \le 
\tfrac12\D_X(\e-\d). \end{equation}
Such numbers $\d$ exist since
$\D_X(\eta)>0$ for all $0<\eta\le 2$.

We argue by contradiction. Suppose there exist $0<\e\le 2$, a number $\d>0$ 
satisfying \eqref{delta},  vectors $x\in S_X$ and
$x'\in X$, and a linear contraction $T$ from $X$ into some Banach space
$Y$ such that the assumptions of (2) are satisfied while
$\n x - x'\n \ge \e.$
From $\d\le \frac12\D_X(\e-\d)\le \frac12$ and {\rm (i)} it follows that
$x'\not=0$. We estimate
\begin{equation}\label{(3)}
\Bigl\n {x} - \frac{x'}{ \n x'\n }\Bigr\n
 \ge {\e}
- \Bigl\n {x'} - \frac{x'}{\n x'\n}\Bigr\n
= \e - \bigl| \n x'\n -1 \bigr| > {\e} - {\d}.
\end{equation}
and 
$$
\tfrac12 \n {x+ x'}\n \le \tfrac12 \Big\n
{x} + \frac{x'}{ \n x'\n }\Big\n + \tfrac12 \Big\n
{x'} -  \frac{x'}{ \n x'\n } \Big\n
 \stackrel{{\rm(a)}}{<} \big(1-\D_X(\e-\d)\big) + \tfrac12{\d}
 \stackrel{{\rm(b)}}{\le}
 1-\tfrac34\D_X(\e-\d).
$$
In (a), the first term is estimated using  \eqref{(3)} and 
the definition of the modulus 
of convexity and the second term is estimated using assumption (i). 
The estimate (b) is immediate from \eqref{delta}.
Thus, $\xi:= x+ x'$ satisfies $\n \xi\n < 2-\frac32\D_X(\e-\d).$ 
Next, from $\xi':= 2Tx- T\xi = T x-Tx'$
and {\rm (iii)} it follows that
$\n \xi'\n  <\d.$
Putting things together and using {\rm (ii)}, we obtain
$$
2 - 2{\d}
< \n 2Tx \n  = \n T\xi+\xi'\n\le \n T \xi\n + \n \xi'\n
 \le \n \xi\n+\n \xi'\n < \big(2-\tfrac32\D_X(\e-\d)\big) + \d.
$$
Comparing left- and right hand sides, we have obtained a contradiction
with \eqref{delta}.

(2)$\Rightarrow$(3): \ Trivial.

(3)$\Rightarrow$(1): \ If $X$ is not uniformly convex,
there exist $\e>0$ and sequences
$(x_n)\subseteq S_X$,  $(x_n')\subseteq S_X$
such that $\inf_n \n x_n - x_n'\n\ge \e$ and
$\limn \n x_n+x_n'\n = 2$.
Choose a sequence $(x_n\s)\subseteq S_{X\s}$ such that
$ \lb x_n+x_n', x_n\s\rb = \n x_n + x_n'\n $ for all $n$.
Then $\limn \lb x_n,x_n\s \rb = \limn \lb x_n', x_n\s\rb =1$ and
$\limn \lb x_n - x_n', x_n\s\rb =0$, and therefore by (3) we obtain
$\limn \n x_n - x_n'\n=0$, a contradiction.
\end{proof}

For the proofs of Theorems \ref{thm:main1} and \ref{thm:main3}
we need some quantitative information about the dependence of $\d$ upon $\e$
in the proof of (1)$\Rightarrow$(2)$_0$, where (2)$_0$
is obtained from (2) by replacing {\rm (ii)} by the more restrictive condition
\begin{enumerate}
\item[]
\begin{enumerate}
\item[(ii)$_0$] $\dps \n T x\n = 1$.
\end{enumerate}
\end{enumerate}

\begin{lemma}\label{lem:23}
Let $X$ be uniformly convex and fix an arbitrary $0<\e\le 2$. 
Then the conclusion of {\rm(2)}$_0$ holds if the
assumptions {\rm(i)}, {\rm(ii)}$_0$, {\rm(iii)} are satisfied for 
$\d =  \D_X(\tfrac23{\e}).$
\end{lemma}
\begin{proof}
We first claim that the conclusion of (2)$_0$ holds if {\rm (i)}, {\rm
(ii)}$_0$, {\rm (iii)}
are satisfied for some $\d>0$ such that
\begin{equation}\label{d_0}
\d\le \D_X(\e-\d).
\end{equation}
Arguing by contradiction and proceeding as in the proof of (1)$\Rightarrow$(2)
we first obtain
$$\tfrac12 \n x+y\n < 1 - \tfrac12\D_X(\e-\d)$$ 
and then, with {\rm (ii)}$_0$,
$$ 2 < \big(2-\D_X(\e-\d)\big)+\d.$$
This contradicts the choice of $\d$ and the claim is proved.

It remains to check that \eqref{d_0} holds for
$\d= \D_X(\frac23\e)$. But from $\n x_1+x_2\n \ge 2\n x_1\n-\n
x_1-x_2\n$ we have, for all $0<\eta\le 2$,
$$ \D_X(\eta) \le \inf\big\{1-\tfrac12\n x_1+x_2\n: \ x_1,x_2\in S_X, \ \n
x_1-x_2\n = \eta\big\} \le \tfrac12\eta.$$
Hence if $\d= \D_X(\frac23\e)$, then
$\d \le \tfrac12\cdot \tfrac23\e = \tfrac13\e$ and consequently,
$\d = \D_X(\e - \tfrac13\e)\le
\D_X(\e-\d)$
by the monotonicity of $\D_X$.
\end{proof}

In a similar way one checks that the conclusion of (2) holds 
if {\rm (i)}, {\rm (ii)}, {\rm (iii)} are satisfied for $\d= \tfrac12\D_X(\tfrac45{\e}).$

\begin{proof}[Proof of Theorem \ref{thm:main1}.]
Assume, for a contradiction, that the theorem were false. Then, for some 
$0<\e\le 2$ and some bounded $\e$-separated sequence $(x_n)$,
there exists an $x\in S_X$ such that every
subsequence $(x_{n_j})$ of $(x_n)$ contains two further subsequences
$(x_{n_{j_k}^{(1)}})$ and $(x_{n_{j_k}^{(2)}})$, with
$n_{j_k}^{(1)}\not=n_{j_k}^{(2)}$ for all $k$, satisfying
\begin{equation}\label{*}
 \big\n x-\big(x_{n_{j_k}^{(1)}}-x_{n_{j_k}^{(2)}}\big)\big\n < 1+\d_\e \quad\hbox{for all $k$},
\end{equation}
where $\d_\e:=\D_X(\frac23 \e).$

Choose $x^{*}\in
S_{X\s}$ with $\lb x,x^{*}\rb =1$. 
Since $(x_n)$ is bounded we may pass to a subsequence
$(x_{n_j})$ for which the limit
$\lim_{j\to\infty} \lb x_{n_j}, x\s\rb$ exists.
We extract two further subsequences $(x_{n_{j_k}^{(1)}})$ and
$(x_{n_{j_k}^{(2)}})$ of $(x_{n_j})$ satisfying \eqref{*} and put
$$\xi_{k} := x- \bigl(x_{n_{j_k}^{(1)}}  - x_{n_{j_k}^{(2)}}\bigr).$$
Then $ \n \xi_k\n < 1+ \d_\e$ for all $k$ and 
$$\lim_{k\to\infty} \lb
x-\xi_{k}, x\s\rb =0.$$ 
Hence $\lim_{k\to\infty}  \lb \xi_{k}, x^{*}\rb
= 1$. In particular,
$\n \xi_k\n > 1-\d_\e$ for large $k$.
Thus,
$$ \big|1- \n \xi_k\n\big|< \d_\e\quad \hbox{for large $k$}.$$
By Lemma \ref{lem:23} the conclusion of (2)$_0$ applies with $T := x\s$.
As a result,
 for large $k$ we obtain
$$
 \big\n x_{n_{j_k}^{(1)}}  - x_{n_{j_k}^{(2)}}\big\n
= \n x-\xi_{k}\n  < \e.
$$
But this contradicts the fact that $(x_n)$ is $\e$-separated.
\end{proof}

\begin{proof}[Proof of Theorem \ref{thm:main3}.]
We start from an arbitrary $1$-separated sequence
$(\xi_n)_{n\in\N}\subseteq S_X$; a short and elementary construction of such 
sequences is given in the notes of \cite[Chapter~1]{Di}.

Let $\d_1:= \D_X(\frac23)$. By Ramsey's theorem \cite{Go}, $(\xi_n)$ has a subsequence $(\xi_{n_j})$ such that
either 
$$\n \xi_{n_j}-\xi_{n_k}\n \in [1,1+\tfrac12\d_1]\quad \hbox{for all
$j\not=k$}$$ or
$$\n \xi_{n_j}-\xi_{n_k}\n \in (1+\tfrac12\d_1,2]\quad \hbox{for all
$j\not=k$.}$$
In the second case we are done (take $x_j = \xi_{n_j}$ and recall that
$\xi_{n_j}\in S_X$).
Hence, after relabeling we may assume that 
\begin{equation}\label{relabeling}
\n \xi_{j}-\xi_{k}\n \in [1,1+\tfrac12\d_1]\quad \hbox{for all $j\not=k$}.
\end{equation}

Let $\phi:\N\to (\N\times \N)\setminus \mathbb{D}$ be a bijection, where
 $\mathbb{D} = \{(j,k)\in\N\times\N: \ j=k\}$, and write
$\phi(n) = (\phi_1(n), \phi_2(n))$. 
Put
$$x_0:= \frac{\xi_{\phi_1(0)} - \xi_{\phi_2(0)}}{\n \xi_{\phi_1(0)} -
\xi_{\phi_2(0)}\n}.$$
Suppose next that integers 
$0=:n_0 < \dots<n_{m-1}$ have been chosen subject to
the condition that the vectors
$$x_j:=\frac{y_j }{\n y_j\n}, \quad 0\le j\le m-1$$
where
$ y_j := \xi_{\phi_1(n_j)} - \xi_{\phi_2(n_j)},$
satisfy
$$\n x_j - x_k\n \ge 1+ \tfrac12\d_1\quad\hbox{for all $0\le j< k\le m-1$.}$$ 
By Theorem
\ref{thm:main1}, applied consecutively to the vectors $x=x_j$, $0\le j\le m-1$, there exists an
integer $n_m > n_{m-1}$ such that 
$$
\n x_j -y_m\n \ge 1+\d_1
\quad\hbox{for  all $\ 0 \le j\le m-1$},
$$
where $ y_m := \xi_{\phi_1(n_m)} - \xi_{\phi_2(n_m)}$.
With 
$$x_{m} := \frac{y_m}{\n y_m\n}
$$
we have, for all $0\le j\le m-1$,
$$
\begin{aligned}
 \n x_j - x_m\n & \ge (1+\d_1) - 
 \Big\n y_m-
\frac{y_m}{\n y_m \n} \Big\n 
\\ & = (1+\d_1) - \big| \n y_m\n -1\big|
\ge (1+\d_1) - \tfrac12 \d_1 = 1+\tfrac12\d_1, 
\end{aligned}
$$
where the last inequality follows from \eqref{relabeling}.
Continuing this way we obtain a
sequence $(x_n)_{n\in\N}$ with the desired properties.
\end{proof}

A Banach space $X$ is called
{\em locally uniformly rotund} \cite{De} 
if for all $x\in S_X$ and all sequences $(x_n)\subseteq S_X$
with $\limn \n x+x_n\n_X = 2$ we have $\limn \n x - x_n\n_X = 0$.
A characterization of locally uniformly rotund Banach  spaces analogous to Theorem \ref{thm2} holds; 
the numbers $\d$ in (2) and (3) will now depend on $\e$ and $x$.
As a result, Theorem \ref{thm:main1} remains true with a separation constant
depending on $x$.

\bigskip

{\em Acknowledgment} -- The author thanks Mark Veraar for pointing out reference \cite{KP}.

\end{document}